\documentclass[a4paper,11pt]{amsart}
\usepackage{t1enc,amssymb,epsfig}
\usepackage[english]{babel}

\input xy
\xyoption{all}

\newtheorem{proposition}{Proposition}[section]
\newtheorem{theorem}[proposition]{Theorem}
\newtheorem{lemma}[proposition]{Lemma}
\newtheorem{definition}[proposition]{Definition}
\newtheorem{corollary}[proposition]{Corollary}
\newtheorem{remark}[proposition]{Remark}

\newcommand{\C}{\mathbf{C}}
\newcommand{\Z}{\mathbf{Z}}

\newcommand{\R}{\mathbf{R}}

\begin{document}

\title[Examples with bijective Nash map]{\textbf{Families of higher
    dimensional germs \\ with bijective Nash map}}
\author{\sc Camille Pl{\'e}nat}
\address{Universit{\'e} de Provence, LATP UMR 6632\\
Centre de Math{\'e}matiques et Informatique, 39 rue Joliot-Curie,13453
    Marseille cedex 13, France.} 
\email{plenat@cmi.univ-mrs.fr}
\author{\sc Patrick Popescu-Pampu}
\address{Univ. Paris 7 Denis Diderot, Inst. de
  Maths.-UMR CNRS 7586, {\'e}quipe "G{\'e}om{\'e}trie et dynamique" \\case
  7012, 2, place Jussieu, 75251 Paris cedex 05, France.}
\email{ppopescu@math.jussieu.fr}



\medskip

\thispagestyle{empty}
\begin{abstract}{Let $(X,0)$ be a germ of complex analytic
    normal variety, non-singular outside $0$. An essential divisor
    over $(X,0)$ is a divisorial valuation of the field of meromorphic
    functions on $(X,0)$, whose center on any resolution of the germ is
    an irreducible component of the exceptional locus. The Nash
    map associates to each irreducible component  of the space of
    arcs through $0$ on $X$ the unique essential divisor intersected 
    by the strict transform of the generic arc in the component.
    Nash proved its injectivity and asked if it was bijective. We
    prove that this is the case if there exists a divisorial
    resolution $\pi$ of $(X,0)$ 
    such that its reduced exceptional divisor carries sufficiently
    many $\pi$-ample divisors (in a sense we define). Then we apply
    this criterion 
    to construct an infinite number of families of 3-dimensional
    examples, which 
    are not analytically isomorphic to germs of toric 3-folds (the
    only class of normal 3-fold germs with  
    bijective Nash map known before).}
\end{abstract}

\maketitle

\par\medskip\centerline{\rule{2cm}{0.2mm}}\medskip
\setcounter{section}{0}

\section{Introduction}

Let  $X$ be a reduced complex algebraic  variety.  An \emph{arc}
contained in $X$ is a germ of formal map: 
$$(\C,0) \rightarrow X.$$
If $t$ denotes the local parameter of $\mathbf{C}$ at $0$, notice that
each arc comes equipped with a canonical  parametrization: thought
algebraically, it is a morphism of $\mathbf{C}$-algebras
$\mathcal{O}_{X,0} \rightarrow \mathbf{C}[[t]]$. 

In a preprint written around 1966, published later as \cite{N 95}, Nash defined
the associated \emph{arc space}   $X_{\infty}$ of $X$, whose points
represent the  arcs contained in $X$. By looking at the Taylor
expansions of the functions on $X$ with respect to the parameter $t$
and to their truncations at all the orders, Nash constructed this
space as a projective limit of algebraic varieties of finite type over
$X$. 

If one associates to a formal arc the point
of $X$ where it is based, that is the image of $0 \in \mathbf{C}$, one
gets a natural map:  
$$\alpha: X_{\infty} \rightarrow X$$
If $Y$ is a closed subvariety of $X$, denote by:
$$(X,Y)_{\infty}:=\alpha^{-1}(Y)$$
the \emph{space of arcs on}  $X$ \emph{based at} $Y$. 

Nash was thinking of the spaces $X_{\infty}$ and $(X,Y)_{\infty}$ for
varying $Y \subset \mathrm{Sing}(X)$ as
tools for 
studying the structure of $X$ in the neighborhood of its singular
set.  Indeed, the main object of his paper was to state a program  for
comparing the various resolutions of the singularities of $X$. Such
resolutions always exist, as had recently been proven by Hironaka, but
unlike in the case of surfaces, minimal ones do not necessarily
exist. We quote from the introduction of \cite{N 95} the two main
problems formulated by Nash in this direction: 

\medskip

i) \emph{For surfaces it seems possible that there are exactly as many
   families of arcs associated  with a point as there are components
   of the image of the point in the minimal  resolution of the
   singularities of the surface}.  

ii) \emph{In higher dimensions, the arc families associated with the
  singular set correspond to "essential components"   which must
  appear in the image of the singular set in all resolutions.  We do
  not know how complete is the representation of essential  components
  by arc families.} 
\medskip

The first question is a local one, as it deals with the structure of
$X$ only in a neighborhood of one of its (closed) points.  The second
one is more global, as it deals with the structure of $X$  in the
neighborhood of its entire singular set.  

Following Nash's paper, the foundations for his program were worked
with more detail by Lejeune-Jalabert \cite{LJ 90}, Nobile \cite{N 91}
and Ishii $\&$ Koll{\'a}r \cite{IK 03}. They also extended the program to
other categories of spaces. For example, Ishii $\&$ K{\'o}llar \cite{IK
  03} considered  
schemes over arbitrary fields, Lejeune-Jalabert 
\cite{LJ 90} and Nobile \cite{N 91} considered  formal germs of
varieties. Their treatment extends readily to germs of complex
analytic varieties.

For such germs, the space $(X, \mathrm{Sing}(X))_{\infty}$ of arcs
based at the singular 
locus of $X$ can be canonically given the
structure of a relative scheme over $X$, as the projective limit of
relative schemes of finite type obtained by truncating arcs at each finite
order. 

\medskip
\emph{In the sequel we will restrict to the case where $(X,0)$ is a germ of
a complex analytic variety and $\mathrm{Sing}(X)= \{0\}$.}

The space $(X, 0)_{\infty}$ of
arcs on $X$ based at $0$ is a relative subscheme over $X$ of $X_{\infty}$. As
it projects onto $0$, we see that it is in fact a true scheme (but not of
finite type over $\mathbf{C}$). 
This implies that it makes sense to speak about the set
$\mathcal{C}(X,0)_{\infty}$ of its irreducible components.

Denote by
$$\pi: \tilde{X} \rightarrow X$$
a resolution of $X$. The \emph{exceptional set} $\mathrm{Exc}(\pi):=
\pi^{-1}(0)$ is not 
assumed to be of pure codimension 1, that is, the resolution is not
necessarily divisorial. 

An irreducible component of $\mathrm{Exc}(\pi)$ is called
\emph{an essential component of $\pi$} if it
corresponds to an irreducible 
component of the exceptional set of any other resolution of $X$.  In 
other words, if its birational transform is an irreducible component of the
exceptional set in any resolution. An equivalence class of such
essential components over all the resolutions of $X$
is called an 
\emph{essential divisor over} $(X,0)$. If we denote by
$\mathcal{E}(X,0)$ the set of 
essential divisors over $(X,0)$, the  essential
components of the given resolution morphism $\pi$ are
in a canonical 
bijective correspondence with the elements of $\mathcal{E}(X,0)$.

Let $\mathcal{K}$ be an element of $\mathcal{C}(X,0)_{\infty}$. For each arc
  represented by a point of $\mathcal{K}$, one can consider the
  intersection point with $\mathrm{Exc}(\pi)$ of its 
strict transform on $\tilde{X}$. For an arc generic with respect to
  the Zariski topology of  $\mathcal{K}$, this intersection 
  point is situated on a unique irreducible component of
  $\mathrm{Exc}(\pi)$; moreover, this component is essential (Nash
  \cite{N 95}). In
  this manner one defines a map:
    $$ \mathcal{N}_{X,0} :  \mathcal{C}(X,0)_{\infty} \rightarrow
    \mathcal{E}(X,0)$$ 
which is called \textit{the Nash map} associated to 
    $(X,0)$. Nash proved that the map
    $\mathcal{N}_{X,0}$ is always \emph{injective} (which shows in
    particular that $\mathcal{C}(X,0)_{\infty}$ is 
    a finite set). In our context, one can reformulate question ii) above:
$$
  {\mbox{\emph{When is the map} } \mathcal{N}_{X, 0}
    \mbox{ \emph{bijective?}} }
$$   

This question is also known as \emph{the Nash problem on arcs}. 

In \cite{PPP 05} we listed the classes of isolated \textit{surface}
singularities for 
which the Nash map was proved to be bijective. In higher dimensions, the
bijectivity of $\mathcal{N}$ was  proved till now for the following 
classes of germs with isolated singularities:

$\bullet$ for the germs which have resolutions with irreducible
exceptional set, for trivial reasons; 

$\bullet$ for germs of normal toric varieties by Ishii and
Koll{\'a}r in \cite{IK 03} (also in the case of non-isolated singularities);

$\bullet$ for quasi-ordinary singularities by Ishii in \cite{I 04} (in
fact for a class slightly more general, which contains also non-isolated
singularities).
\medskip

No surface  or 3-fold is known for which the Nash map is not bijective. 
But Ishii and Koll{\'a}r
proved in \cite{IK 03}  that it is not always bijective for algebraic
varieties of dimension at least 4.  Indeed, they gave a counterexample
in dimension $4$, which can be immediately 
transformed into a counterexample (with non-isolated singularity) in
any larger dimension. 

\medskip

In this article we construct \emph{a class of normal 
isolated singularities of arbitrary dimension} $(X,0)$ for which the Nash map
$\mathcal{N}_{X,0}$ is bijective (Corollary
\ref{amplebij}). The definition of the class uses  a 
criterion ensuring that a divisorial component of the exceptional set
of a given resolution is in the image of the Nash map (Theorem
\ref{maincrit}). In fact, we use that theorem in the following less
general form (a reformulation of Corollary \ref{gsam}), which allows
us to apply Kleiman's ampleness criterion:

\medskip

\textbf{Theorem} \textit{Let $\pi:\tilde{X}\rightarrow X$ be a
  divisorial projective resolution of $(X,0)$. Consider an irreducible
  component $E_i$ of $\mathrm{Exc}(\pi)$. Suppose that for any other
  component $E_j$, there exists an effective integral divisor $F_{ij}$
  on $\tilde{X}$ 
  whose support coincides with $\mathrm{Exc}(\pi)$, in which the
  coefficient of $E_i$ is strictly less than the coefficient of
  $E_j$ and such that the line bundle
  $\mathcal{O}_{\mathrm{Exc}(\pi)}(-F_{ij})$ is ample. Then $E_i$ is
  an essential component contained in the image 
  of the Nash map.}

\medskip

Using the previous criterion, we construct an infinite
family of examples of 3-dimensional singularities with bijective Nash map (see
Section \ref{examples}). In Section  \ref{aninv}, we distinguish some
of the singularities constructed before using suitable analytical
invariants. Moreover, we determine those which are isomorphic to germs
of toric varieties, establishing like this the intersection of our
class of examples with the classes known before.

\subsection*{Acknowledgements} We are grateful to Monique
Lejeune-Jalabert and to the referee of \cite{PPP 05} for having
remarked that the main theorem of \cite{PPP 05} could be obtained by
replacing Laufer's vanishing theorem with a convenient use of the
ampleness of adequate line bundles. This is the method we have
extended in this paper to higher dimensions. The second author is also
grateful to S{\'e}bastien Boucksom, David Eisenbud, Charles Favre and
Christophe Mourougane for their explanations about ampleness and
positivity. He addresses special thanks to C. Mourougane for his many
helpful remarks on previous versions of this paper. 

\bigskip

\section{Essential divisors and essential components} \label{critfund}

\medskip

In the sequel, if $A$ is a complex analytic space or a relative scheme
over an analytic space, we denote by
$\mathcal{C}A$ the set of its irreducible components. 

Let $(X,0)$ be an irreducible  germ of  complex analytic
variety. We suppose 
that $\mathrm{Sing}(X)=\{0\}$, that is, the germ is smooth outside the origin
(with a slight abuse of vocabulary due to the fact that $X$ is also allowed
to be smooth, we say that the germ \emph{has an isolated
  singularity}).  Denote by 
$\frak{m}$ the maximal ideal of its local ring
$\mathcal{O}_{X,0}$. We also write $\mathcal{N}$ instead of
$\mathcal{N}_{X,0}$, as we do not  consider various Nash maps at the
same time. 

Consider a \emph{resolution}  $\pi: \tilde{X} 
\rightarrow X$. This means that $\pi$ is a proper bimeromorphic
map with $\tilde{X}$ smooth, restricting to an isomorphism over the
complement of $0$ in $X$. The \emph{exceptional set
  $\mathrm{Exc}(\pi)$ of $\pi$} is by 
definition the subset of $\tilde{X}$ where  $\pi$ is not a local
isomorphism. If $0$ is a singular point of $X$, it coincides with the
preimage $\pi^{-1}(0)$. If each irreducible component of  
$\mathrm{Exc}(\pi)$ is of pure codimension 1 in $\tilde{X}$, we say
that $\pi$ is \emph{divisorial}. In the sequel, we do not suppose
that this is the case. We do neither suppose that the
morphism $\pi$ is projective. 

\begin{remark} \label{smallres}
 In dimension $2$, all the
resolutions of a normal surface are divisorial. This is no longer true
in higher dimensions: the simplest example of a normal germ with
isolated singularity which has non-divisorial resolutions is the
3-fold hypersurface germ defined by the affine equation $xy-zt=0$.
Nevertheless,  all the 
resolutions of a $\mathbf{Q}$-factorial germ are divisorial (see
Debarre \cite[Section 1.40]{D 01}).   
\end{remark}

Consider a closed irreducible subvariety $E$ of $\mathrm{Exc}(\pi)$ (not
necessarily one of its irreducible components). Take  
the preimage $D$ of $E$ on $B_E(\tilde{X})$, the 
variety obtained by blowing-up $E$ in $\tilde{X}$. As
$\tilde{X}$ is smooth, this preimage is an irreducible hypersurface of
$B_E(\tilde{X})$. Therefore, it induces a discrete  valuation
$v_E$ of
rank 1  on the field of meromorphic functions on $(X,0)$  (which
associates to any such function the order of vanishing along $D$ of
its total transform on $B_E(\tilde{X})$). 

If $\psi: \overline{X}\rightarrow X$ is another resolution of $X$,
\emph{the birational transform} $E^{\psi}$ of $E$ on $\overline{X}$ is the
center of the valuation $v_E$ on $\overline{X}$.  
We have obviously $v_E = v_{E^{\psi}}$. This allows to
identify the valuation $v_E$ with the set whose elements are $E$ and
its birational transforms on all the resolutions of $X$. Following
\cite{IK 03}, we say
that $v_E$ (or the class of its centers on all the resolutions) is an 
\emph{exceptional divisor over $(X,0)$}. The name is motivated by the
fact that any 
resolution is dominated by another one on which the center of $v_E$ is
a divisor (as above, just blow-up $E$, then resolve the singularities
of the new space). 

Conversely, if $v$ is an exceptional divisor over $(X,0)$ and $\pi:
\tilde{X} \rightarrow X$ is a resolution, we denote by $E_v^{\pi}$ (or
$E_v$, if $\pi$ is clear from the context) the center of $v$ on
$\tilde{X}$. Among the exceptional divisors over $(X,0)$, Nash  
distinguished those whose centers are not only subvarieties, but
irreducible 
\emph{components} of the exceptional locus of any resolution of
$(X,0)$ (in fact he considered this in the global case of an algebraic
variety; Ishii \cite[Definition 2.10]{I 04} considers the same
localized situation as ours):

\begin{definition} \label{essdiv}
An \textbf{essential divisor over} $(X,0)$ is an exceptional divisor
over $(X,0)$ whose center on $\tilde{X}$ is an irreducible component
of $\mathrm{Exc}(\pi)$, this for any resolution $\pi:\tilde{X} \rightarrow
X$. We also say that the centers of the essential divisors on
$\tilde{X}$ are the \textbf{essential components} of $\pi$. 
\end{definition}

In \cite{B 98}, Bouvier considered another definition of
\emph{essential divisors}. She called a component of codimension 1 of
the exceptional set \emph{essential} if its center on any resolution
was a divisor. Her definition is strictly more restrictive than ours,
as shown by the germs which admit resolutions without exceptional
components of codimension 1 (see Remark \ref{smallres}). Ishii and
Koll{\'a}r introduced a third notion in \cite{IK 03}, that of
\emph{divisorially essential divisors}. Namely, an exceptional divisor
is of this type, if its center in any divisorial resolution is an
irreducible component of the exceptional set. It follows directly from
the Definition \ref{essdiv} that an essential divisor is a
divisorially essential divisor, but it seems to be an open question if
the converse is true.

\medskip
\emph{In the sequel we consider only the notion of essential components and
essential divisors introduced in Definition \ref{essdiv}}.

If $(X,0)$ is a normal surface singularity, then the essential
divisors over $(X,0)$ are precisely the divisorial valuations
generated by the irreducible components of the exceptional set of the
minimal resolution of $(X,0)$. In higher dimensions it is much more
difficult to determine them, as no minimal resolution (in the sense
that it is dominated by all the other resolutions) exists in general.

The only class of singularities for which the essential divisors are
completely known is that of germs of normal toric varieties. Indeed, Bouvier
\cite{B 98} determined combinatorially the 
essential divisors of all normal toric germs. Her work was based on
preliminary results of Bouvier $\&$ Gonz{\'a}lez-Sprinberg \cite{BGS 95}. 

Two general criteria are known,  ensuring that a 1-codimensional
component of the exceptional locus of a given resolution is essential
(see Ishii $\&$ Koll{\'a}r \cite[Examples 2.4, 2.5, 2.6]{IK 03}):

\begin{proposition}
  Let $E_i$ be an irreducible component of $\mathrm{Exc}(\pi)$, which
  is of codimension 1 in $\tilde{X}$. 

 1) (Nash \cite{N 95}) If $E_i$ is not birationally
  ruled, then $E_i$ is essential.

 2) If $(X,0)$ is a canonical singularity and $E_i$ is crepant, then
    $E_i$ is essential.  

Moreover, in both cases the birational transform of $E_i$ on any other
resolution has again codimension 1. 
\end{proposition}

One of
the results of our work is to give a new  criterion of essentiality
for exceptional divisors, using
the space of arcs on $X$ based at $0$ (Theorem
\ref{maincrit}).

\medskip

For each irreducible component $E$ of $\mathrm{Exc}(\pi)$, consider the smooth 
arcs on $\tilde{X}$ 
whose closed points are on $E- \cup_{F \neq E}F$, where $F$ varies
among the elements of $\mathcal{C}\mathrm{Exc}(\pi)$, and which
intersect $E$ transversely (that is, such that their tangent line and
the tangent space to $E$ at their intersection point are direct
summands). Consider the set of their images in 
$(X,0)_{\infty}$ and denote the closure of this set by
$V(E)$. 

\begin{remark} 
  In fact $V(E)$ only depends on the exceptional divisor $v_E$ over
  $(X,0)$ determined by $E$ (see Ishii \cite[Example 2.14]{I 04}). For this
  reason, in the sequel we also write $V(v)$ instead of
  $V(E)$, if $E=E_v$. 
\end{remark}

Nash \cite{N 95} proved:

\begin{proposition} \label{decomp}
  1) The sets $V(E)$  are irreducible subvarieties
     of $(X,0)_{\infty}$ (but not necessarily components).

  2) $$(X,0)_\infty = \bigcup_{E\in \mathcal{C}\mathrm{Exc}(\pi)}V(E).$$
\end{proposition}

\medskip

The next lemma gives a  criterion to show that an
exceptional divisor $v$ over $(X,0)$ is essential, using its image
$V(v)$ in the space of arcs based at $0$.

\begin{lemma} \label{elemess}
 Let $v$ be an exceptional divisor over $0$. If $V(v)$ is an
 irreducible component of $(X,0)_{\infty}$, then $v$ is 
 essential.
\end{lemma}
\begin{proof}
  Suppose by contradiction that $v$ is inessential. This means that
  there exists a resolution $\pi:\tilde{X}\rightarrow X$ such that the
  center $E_v$ of the valuation $v$ on $\tilde{X}$ is strictly
  included in an irreducible component $E$ of $\mathrm{Exc}(\pi)$. We deduce
  that $V(v)$ is strictly included in $V(E)$. But this
  last variety is irreducible, by Proposition \ref{decomp}. This
  contradicts the fact that $V(v)$ is an irreducible component
  of $(X,0)_{\infty}$. 
\end{proof}

The next proposition gives a criterion to prove that some
components of the exceptional locus of a resolution of $(X,0)$ are
essential, and in particular to prove that Nash's map
$\mathcal{N}$ is bijective.

\begin{proposition} \label{critess}
  Let $\pi:\tilde{X}\rightarrow X$ be a resolution. Consider the
  irreducible components $(E_i)_{i\in I}$ of $\mathrm{Exc}(\pi)$. Suppose
  that one can write the index set $I$ as a disjoint union $I = J
  \bigsqcup K$ such 
  that $V(E_j) \not \subset  V(E_i), \: \forall \: j \in
  J, \: \forall \: i \in I-\{j\}$. Then:

   1) The varieties $(V(E_j))_{j \in J}$ are
   irreducible components of $(X,0)_{\infty}$. In particular,
   $(E_j)_{j \in J}$ are essential components of $\pi$. 

   2) If  $(E_k)_{k \in K}$ are all inessential components of $\pi$, then
  $\mathcal{N}$ is bijective. 
\end{proposition}

\textbf{Proof}
 1) By Proposition \ref{decomp},  the irreducible components of
   $(X,0)_\infty$ are among the varieties $V(E)$ with $E\in
   \mathcal{C}E_{\pi}(X,0)$. Moreover, by definition, the irreducible
   components of $(X,0)_\infty$ are those which are not included in other
   irreducible subsets. Then the varieties $(V(E_j))_{j \in J}$
   are 
   irreducible components of $(X,0)_{\infty}$. By the previous lemma,
   it follows that $(E_j)_{j \in J}$ are essential components of $\pi$. 

 2) If the components $(E_k)_{k \in K}$ are all inessential, then,
  by Lemma \ref{elemess}, the varieties  $(V(E_k))_{k \in K}$ are not
  irreducible components of  $(X,0)_\infty$. The irreducible
  components of $(X,0)_\infty$ are exactly $(V(E_j))_{j \in J}$
  and they do correspond bijectively to the essential components of
  the resolution. This implies that $\mathcal{N}$ is bijective.  
   \hfill $\Box$

\medskip

The following proposition was proven  
by the first author in \cite[2.2]{P 03}, as a generalization of
Reguera \cite[Theorem 1.10]{R 95}, who considered only the class of
rational surface singularities. It is an essential ingredient
of all the  criteria we prove in this paper. It was also the
basis of our work \cite{PPP 05}.  

\begin{proposition} \label{distinguish}
    Let $v_1$ and $v_2$ be exceptional divisors over $(X,0)$. If there
   exists a function $f \in \frak{m}$ such 
   that  $v_1(f) < v_2(f)$, then $V(v_1) \nsubseteq
   V(v_2)$.
\end{proposition}

In Section \ref{criterion}, we combine the propositions \ref{critess} and
\ref{distinguish} in order to give criteria of
essentiality for exceptional divisors in terms of global generation
and ampleness of suitable line bundles. Before that, we need some
background about ampleness and exceptional sets.

\bigskip

\section{Background about ampleness and exceptional analytic sets}
\label{contractions}

In this section we recall Kleiman's criterion of ampleness and
Grauert's criterion of contractibility. 

Let $Y$ be a complete algebraic variety. Let $Z_1(Y)_\R $  be the $\R $-vector
space of real one cycles on X, consisting of all finite  $\R$-linear
combination of irreducible algebraic curves on $Y$. Two elements
$\gamma_1$ and $\gamma_2$ of $Z_1(Y)_\R$ are \emph{numerically
  equivalent} if one has the equality of intersection numbers
$$E\cdot \gamma_1\: =\: E\cdot \gamma_2$$ for every $E\in 
\mathrm{Div}(Y)\otimes_{\Z} \R $, where $\mathrm{Div}(Y)$ denotes the
group of Cartier divisors on $Y$.  The 
corresponding vector space of numerical equivalence classes  of
one-cycles is written $N_1(Y)_\R$. 

\begin{definition}\label{cocu}
Let $Y$ be a complete algebraic variety. The \textbf{cone of curves}
$$NE(Y)\subset 
N_1(Y)_\R $$ is the cone  $\R_+-$spanned by the classes of all effective
one-cycles on $Y$. 
  Its closure $\overline{NE}(Y)\subset N_1(Y)_\R $ is the \textbf{ closed
  cone of curves} or \textbf{Kleiman-Mori cone} of $Y$.  
\end{definition}

\begin{theorem}{\bf (Kleiman's criterion of ampleness)}\label{klei}
Let $Y$ be  projective variety. 
A Cartier divisor $E$ on $Y$ is ample if and only if $E\cdot z>0$ for all
non zero $z\in \overline{NE}(X)$. 
\end{theorem}

For details, we refer to  
Debarre \cite{D 01} and Lazarsfeld \cite{L 04}.

Ampleness on a reducible variety can be
tested on its irreducible components (see for example Lazarsfeld
\cite[proposition 1.2.16]{L 04}):  
 
\begin{proposition} \label{amplerestr}
  Let $Y$ be a projective variety and $L$ a line bundle on
  $Y$. Then $L$ is ample on $Y$ if and only if the restriction of $L$
  to each irreducible component of $Y$ is ample.
\end{proposition}

We took the following definition from Peternell \cite[definition 2.8]{P 94}:

\begin{definition} \label{except}
  Let $Y$ be a reduced complex space and $E\subset Y$ a compact
  nowhere discrete and nowhere dense analytic set. $E$ is called
  \textbf{exceptional} (in $X$) if there is a complex space $Z$ and a proper
  surjective holomorphic map $\phi:Y \rightarrow Z$ such that:

(1) $\phi(E)$ is a finite set;

(2) $\phi: Y \setminus E \rightarrow Z \setminus \phi(E)$ is biholomorphic;

(3) $\phi_*(\mathcal{O}_Y)=\mathcal{O}_Z$.

Then one says that $\phi$ \textbf{contracts} $E$ (in $Y$).
\end{definition}

One can show that a map $\phi$ which contracts $E$ in $Y$ is \emph{unique} in
the following sense: if $\phi_k : Y \rightarrow Z_k, \: k \in \{1,2\}$
both contract $E$ in $Y$, then there exists a unique analytic
isomorphism $u :Z_1 \rightarrow Z_2$ such that $\phi_2 = u \circ
\phi_1$. 

In the minimal model theory of algebraic
varieties, one considers more general contractions, which are not
necessarily birational maps.  

The vocabulary is coherent with the one used in section
\ref{critfund}. Indeed, if $\pi:\tilde{X} \rightarrow X$ is a
resolution of $X$, then its \emph{exceptional set} $\mathrm{Exc}(\pi)$ in the
sense of section \ref{critfund} is \emph{exceptional in} $\tilde{X}$
in the sense of Definition  \ref{except}, and $\pi$ contracts 
$\mathrm{Exc}(\pi)$ in $\tilde{X}$.

The strategy we use for constructing examples of 3-dimensional
singularities with bijective Nash maps (see Section \ref{examples}) 
works thanks to Grauert's  fundamental criterion of
contractibility (Grauert \cite{G 62}, see
also Peternell \cite[theorem 2.12]{P 94}). A particular case of it is
sufficient  for our purposes:

\begin{theorem} \label{amplecor} {\bf(Grauert's criterion of
    contractibility)} 
  Let $Y$ be a complex manifold and let $E$ be a reduced projective (not
  necessarily smooth or irreducible) hypersurface in $Y$. Suppose that
  there exists an effective divisor $A$ whose support is $E$, such
  that the restriction $\mathcal{O}_E(-A)$ of the line bundle
  $\mathcal{O}_Y(-A)$ to $E$ 
  is ample. Then the analytic hypersurface $E$ is exceptional in $Y$.
\end{theorem}

\begin{remark} 
 1) If $Y$ is a surface, the converse of the theorem is also true. In
    this case, the hypothesis about the existence of $A$ is equivalent
    to the fact that the intersection matrix of $E$ is negative
    definite. For surfaces, the hypothesis of Grauert's criterion of
    contractibility is usually expressed in this last form.  

 2)  The converse of Theorem \ref{amplecor} is not true in a naive form if
  $\mathrm{dim}_{\mathbf{C}} Y \geq 3$, as shown
  by examples of Laufer \cite{L 81} (see also Peternell \cite[Example
  2.14]{P 94}). Nevertheless, there exists a converse if one replaces
  the search of an ample line bundle by that of a coherent sheaf
  $\mathcal{I}$ such that $\mathrm{supp}(\mathcal{O}_Y/ \mathcal{I})=
  E$ and $\mathcal{I}/\mathcal{I}^2$ is positive (see Peternell
  \cite[Theorem  2.15]{P 94}). 
\end{remark}

\bigskip

\section{Criteria for an exceptional divisor to be essential}
\label{criterion} 

Recall  that $(X,0)$ is supposed to be an irreducible germ with
isolated singularity. \emph{From now on, we suppose moreover that $(X,0)$ is 
  normal}. We need this condition in order to be able to conclude that
a bounded holomorphic function on $X\setminus 0$ extends to a function
holomorphic over $X$.

Let $\pi: (\tilde{X}, \mathrm{Exc}(\pi)) \rightarrow
(X,0)$ be a \textit{divisorial} resolution of
$(X,0)$.  Denote by $(E_i)_{i \in I}$ the  irreducible components of
$\mathrm{Exc}(\pi)$.  

Let 
 $$L(\pi):= \bigoplus_{i \in I} \mathbf{Z} E_i$$
be the lattice freely generated by the $E_i$, that is, the lattice of
divisors on $\tilde{X}$ supported by $\mathrm{Exc}(\pi)$. Inside the
associated  
real vector space $L_{\mathbf{R}}(\pi)$, consider the
closed regular cone:
 $$\sigma(\pi):= \bigoplus_{i \in I}\mathbf{R}_+ E_i$$
of the effective $\mathbf{R}$-divisors on
$\tilde{X}$ supported by $\mathrm{Exc}(\pi)$.

For each pair $(i,j) \in I^2$ with $i \neq j$, consider the closed
convex sub-cone $\sigma_{ij}(\pi)$ of $\sigma(\pi)$ defined by:

  $$ \sigma_{ij}(\pi):= \{\sum_{k \in I} a_k E_k \in \sigma(\pi) \: |
  \: a_i \leq 
  a_j\}$$

\begin{theorem} \label{maincrit}
 Fix $i \in I$. 
  Suppose that for each $j \in I \setminus \{i\}$, the cone
  $\sigma_{ij}(\pi)$  contains in its interior an integral
  divisor $F_{ij}$ such that $\mathcal{O}_{\tilde{X}}(-F_{ij})$ is
  generated by its global sections. Then $V(E_i)$ is in the
  image of the Nash map  
  $\mathcal{N}$. In particular, $E_i$ is an essential component
  of $\pi$.  
\end{theorem}

\textbf{Proof:} Consider $F_{ij}\in \mathrm{int}(\sigma_{ij}(\pi))$ such that
  $\mathcal{O}_{\tilde{X}}(-F_{ij})$ is   generated by its global
  sections. Let us consider for a moment
  $\mathcal{O}_{\tilde{X}}(-F_{ij})$ not as a line bundle, but as the
  subsheaf of the structure sheaf $\mathcal{O}_{\tilde{X}}$ formed by the
  holomorphic functions vanishing along $\mathrm{Exc}(\pi)$ at least
  as much as indicated by the coefficients of $F_{ij}$.   

  If $\mathcal{O}_{\tilde{X}}(-F_{ij})$ is generated by global sections,
  then there exists a function $f_{ij} \in H^0(\tilde{X},
  \mathcal{O}_{\tilde{X}}(-F_{ij}))$ whose divisor has a compact part 
  coinciding with $F_{ij}$. 

  As $\pi$ realizes an isomorphism between $\tilde{X} \setminus
  \mathrm{Exc}(\pi)$ and 
  $X \setminus 0$, there exists a function $g_{ij}$ on $X$ vanishing at $0$, 
  continuous on $X$,  
  holomorphic on $X \setminus 0$ and such that $f_{ij}=
  \pi^*(g_{ij})$. As $X$ is supposed to be normal at $0$ (see the
  beginning of the section), we 
  deduce that $g_{ij} \in  \frak{m}$.

 By construction, $v_{E_i}(g_{ij}) < v_{E_j}(g_{ij})$. Proposition
 \ref{distinguish} implies then that $V(E_i) \nsubseteq
   V(E_j)$. As this is true for any pair $(i,j)\in I^2$ with $i \neq
 j$, the proposition follows from Proposition \ref{critess}. \hfill $\Box$ 

\medskip

The following corollary is a direct consequence of the theorem. We
state it as a separate result, in order to be able to use Kleiman's
criterion of ampleness in combination with it. 

\begin{corollary} \label{gsam}
 Fix $i \in I$. 
  Suppose that for each $j \in I \setminus \{i\}$, the cone
  $\sigma_{ij}(\pi)$ contains an integral
  divisor $F_{ij}$ such that $\mathcal{O}_{\tilde{X}}(-F_{ij})$ is
  ample when restricted to each component of $\mathrm{Exc}(\pi)$. Then
  $V(E_i)$ is in  the image of 
  $\mathcal{N}$ and $E_i$ is an essential component of $\pi$
  relative to $0$. 
\end{corollary}

\textbf{Proof:} As  ampleness is an open condition with respect to the
  topology of $L(\pi)$, we see that the hypothesis implies that there
  exists an $F_{ij}\in \mathrm{int}(\sigma_{ij}(\pi))$ such that
  $\mathcal{O}_{\tilde{X}}(-F_{ij})$ is 
  ample when restricted to each component of $\mathrm{Exc}(\pi)$. By
  Proposition  \ref{amplerestr}, it is also ample 
  when restricted to $\mathrm{Exc}(\pi)$. This implies 
  that $\mathcal{O}_{\tilde{X}}(-F_{ij})$ is ample on a neighborhood of
  $E$ in $\tilde{X}$. But then there exists a multiple $-n_{ij} F_{ij}$ of
  the divisor $-F_{ij}$ (where $n_{ij} \in \Z_{>0}$) whose associated
  sheaf is very ample, which implies that
  $\mathcal{O}_{\tilde{X}}(-n_{ij}F_{ij})$ is generated by its global
  sections. 

  The divisor $n_{ij} F_{ij}$ is interior to the cone $\sigma_{ij}(\pi)$, as
  $F_{ij}$ was supposed to be so. This implies that the hypothesis of
  Theorem \ref{maincrit} are satisfied. The conclusion follows.
  \hfill $\Box$
\medskip

A second corollary gives the  criterion of bijectivity of
the Nash map announced in the introduction: 

\begin{corollary} \label{amplebij}
  Suppose that for each pair $(i,j) \in I^2$ with $i \neq j$, the cone
  $\sigma_{ij}(\pi)$  contains  an integral
  divisor $F_{ij}$ such that $\mathcal{O}_{\tilde{X}}(-F_{ij})$ is
  ample when restricted to each component of $E$. Then the components
  of $E$ are 
  precisely the essential components over $0$ and the Nash map
  $\mathcal{N}$ is bijective.
\end{corollary}

\textbf{Proof:} This is an immediate consequence of Corollary
  \ref{gsam}. 
  \hfill $\Box$

\begin{remark}
 When $(X,0)$ is a germ of normal surface and $\pi : \tilde{X}
  \rightarrow X$ is a resolution, the set of effective 
  divisors $F \in \sigma(\pi)$ such that the line bundle
  $\mathcal{O}_{\tilde{X}}(-F)$ is ample is precisely what we called
  \emph{the strict Lipman semigroup}  in \cite[Remark 4.4]{PPP 05} 
  (see Lipman \cite[10.4 and  proof 
  of 12.1 (iii)]{L 69}). Then Corollary \ref{amplebij} restricted
  to germs of  normal 
  analytic surfaces gives exactly the class of singularities found in
  \cite{PPP 05}.  As explained in the Acknowledgements, the present
  work grew out from the wish to generalize the results of that
  article to higher dimensions.  
\end{remark}

\bigskip

\section{An infinite number of families of examples in dimension 3}
\label{examples} 

Corollary \ref{amplebij} gives a method to construct examples of
singularities $(X,0)$ for which the Nash map $\mathcal{N}$ is
bijective. Namely, one starts from a divisorial resolution of
a germ such that the components of the exceptional locus have closed
cones of curves of finite type. The condition on an effective divisor
supported by the exceptional set to have an ample opposite in
restriction to the exceptional set translates then into a finite system
of linear inequalities. If this system has solutions inside all the
cones considered in the corollary then, using the 
corollary,  one has an example with bijective Nash map. 

One could try to start from germs defined explicitly by equations and to use
one of the available algorithms of resolution. Nevertheless, those
algorithms do not allow to compute the
closed cone of curves of a component of the exceptional set. 

For this reason we decided to work differently. The strategy we followed
was to start from smooth projective varieties $S_i$  with 
cones of curves which are closed and of finite type. Then  choose line
bundles over the varieties $S_i$  
with ample duals and glue analytically the total spaces of those line
bundles along neighborhoods of suitable hypersurfaces of the
$S_i$. Of course, the first thing to adjust in order to do such a
gluing,  
is to make a pairing of the  
chosen hypersurfaces and to fix isomorphisms between the elements in
each pair. 

If the gluing succeeds, one gets a smooth analytic
variety $X$ which contains a divisor $E$ obtained topologically by
identifying the chosen pairs of hypersurfaces. The choices should be
done in order to make $E$ \emph{exceptional} in $X$, in the sense of
Definition \ref{except}. Then try to
construct the divisors $F_{ij}$ verifying the conditions of  Corollary
\ref{amplebij}. The hypothesis on the finiteness of the cones of
curves ensures, as explained before, that this search amounts to the
resolution of a finite system of inequalities.

In this section we apply this strategy  to construct an infinite
number of families 
of 3-dimensional examples with bijective Nash map. All of them are
defined by contracting (using Grauert's criterion) a divisor with two
components inside a smooth algebraic threefold obtained by gluing
algebraically 
along Zariski-open sets the total spaces of suitable line bundles over
geometrically ruled surfaces. Both surfaces are obtained by
compactifying total spaces of suitable line bundles over the same
irreducible smooth projective curve. After the plumbing, the two
surfaces meet transversely along a curve which is isomorphic to the starting
curve. We emphasize the fact that this starting curve is \emph{any}
irreducible  smooth projective curve. 
\medskip

In the sequel we say that an algebraic surface $S$ is \emph{geometrically
  ruled} over a curve $C$ if $S$ is the total space of an algebraicaly
  locally trivial bundle over $C$, with fibers projective lines. We
  say that $S$ is \emph{birationally ruled} if it is birationally
  equivalent to a geometrically ruled surface.

Let $C$ be a smooth irreducible projective curve. Consider two
algebraic line bundles $L_1, L_2$ over $C$ such that:
\begin{equation} \label{degree} 
  \mbox{deg}_{C}L_i= -d_i, \: \forall \: i \in \{1,2\}
\end{equation}

We suppose moreover that:
\begin{equation} \label{neg} 
  d_i >0 , \: \forall \: i \in \{1,2\}
\end{equation}

Denote by $E_i$ the total space of the line bundle $L_i$ and by $C_i$
the image of the zero section of $L_i$ in $E_i$. The relations
(\ref{degree}) and (\ref{neg}) imply:
\begin{equation} \label{selfint} 
  C_i \cdot_{E_i} C_i =- d_i <0 , \: \forall \: i \in \{1,2\}
\end{equation}
(the notation $\cdot_{Y}$ means that one considers intersection
numbers inside the smooth space $Y$). 

One can compactify $E_i$ by adding a curve $\tilde{C}_i$ at infinity, 
getting like this a smooth projective surface $S_i$, which is
geometrically ruled over $C$. Denote by $\pi_i$ the morphism:
 $$\pi_i : S_i \rightarrow C$$
which extends the fibration morphism from $E_i$ to $C$. 

By (\ref{selfint}), one gets:
  $$\tilde{C}_i \cdot_{S_i} \tilde{C}_i = d_i$$

For $i\in\{1,2\}$, we consider the following line bundle on the
geometrically ruled surface $S_i$: 
\begin{equation} \label{secline}
  H_i := \mathcal{O}_{S_i}(-x_i \tilde{C}_i) 
       \otimes_{\mathcal{O}_{S_i}} \pi_i^*(L_j)
\end{equation}
where $\{i,j\}=\{1,2\}$. 

It is important to notice that, as an ingredient of
the construction, we pull back one line bundle over $C$ to the total space of a
compactification of the second line bundle. The important thing is that
the total spaces of the restricted line bundles $\pi_1^*(L_2) |_{ E_1}$
and $\pi_2^*(L_1) | _{E_2}$ are canonically isomorphic (see below the
explanation of relation (\ref{canisom})), which allows to glue 
them. But $\pi_i^*(L_j)$ has not an ample inverse on $S_i$, as its
degree on a fiber of the ruling is $0$. This obliges us to twist the
line bundle. We want to do this without changing the crucial
property of the isomorphism of the total spaces of the rectrictions to
$E_i$. That is why we twist with a line bundle having a section whose
divisor is supported by the curve at infinity. 

We pass now to the needed computations. In the definition
(\ref{secline}), the integer $x_i$ is chosen such 
that the following condition is satisfied: 
\begin{equation} \label{cond}
  \check{H}_i \: \mbox{is ample on} \: S_i, \: \forall \: i \in \{1,2\}
\end{equation} 

As $C_i \cdot_{S_i} C_i <0$ (see relation (\ref{selfint})), one
has (see Debarre \cite[1.35]{D 01}):
 $$\overline{NE}(S_i)=NE(S_i)= 
       \R_+[C_i]\oplus \R_+[F_i]$$
where $[F_i]$ is the class of the fibers of the ruling $\pi_i$. By
       Kleiman's  criterion of ampleness \ref{klei}, condition (\ref{cond}) is
       equivalent to:
$$ \left\{ \begin{array}{l}
              \mbox{deg}_{C_i}H_i <0\\
              \mbox{deg}_{F_i}H_i <0
           \end{array}\right.$$
But:
$$ \mbox{deg}_{C_i}H_i = -x_i \tilde{C}_i \cdot_{S_i} C_i + 
       \mbox{deg}_{C_i}\pi_i^*(L_j)= 0 + \mbox{deg}_{C}L_j
      =  -d_j$$
We have used the fact that the curves $C_i$ and $\tilde{C}_i$ are disjoint,
       the projection formula and relation
       (\ref{degree}). In the same manner, using the projection formula and
       the fact that the curves $\tilde{C}_i$ and $F_i$ meet transversely at 
       one point of $S_i$, we get:
$$\mbox{deg}_{F_i}H_i = -x_i$$

As $d_j >0$ by the hypothesis (\ref{neg}), we see that the condition
(\ref{cond}) is equivalent to:
\begin{equation} \label{sign}
     x_i >0, \: \forall \: i \in \{1,2\}
\end{equation}

The line bundle $\mathcal{O}_{S_i}(-x_i \tilde{C}_i)$ is equiped by
construction with a meromorphic section $s_i$ whose divisor is exactly
$-x_i \tilde{C}_i$. This implies that the restriction of $s_i$ to
$E_i=S_i \setminus \tilde{C}_i$ is  a regular and  never 
vanishing section of the restricted line bundle
$\mathcal{O}_{S_i}(-x_i \tilde{C}_i)|_{E_i}$. We deduce that this last
line bundle is trivial. As a consequence:
\begin{equation} \label{structrestr}
   H_i|_{E_i} \stackrel{(\ref{secline})}{=} 
   \mathcal{O}_{S_i}(-x_i \tilde{C}_i)|_{E_i} 
       \otimes_{\mathcal{O}_{E_i}} \pi_i^*(L_j)
   \simeq \pi_i^*(L_j)|_{E_i}
\end{equation}

Denote by $M_i$ the total space of the line bundle $H_i$ over
$S_i$ and by $N_i$ the total space of the line bundle $H_i$ over
$E_i$. Consequently, $N_i$ is a Zariski open set of $M_i$. By relation
(\ref{structrestr}), $N_i$ is isomorphic to the total space of the
line bundle $\pi_i^*(L_j)$ over $E_i$, which in turn is isomorphic to
the total space of the split vector bundle $L_1 \oplus L_2$
of rank 2 over $C$. This gives canonical isomorphisms:
\begin{equation}\label{canisom}
   N_1 \simeq E \simeq N_2
\end{equation}

If we glue the algebraic manifolds $M_1$ and $M_2$ by
identifying $N_1$ and $N_2$, we obtain a new 3-dimensional
algebraic manifold $M: = M_1 \cup M_2$ (with a slight abuse of
notations), in which $S_1$ and $S_2$ are canonically
embedded. We will consequently use the same notation for their images
in $M$. Then:
  $$S_1 \cap S_2= C$$ 
where $C$ is the curve obtained by the identification under the
preceding gluing of the curves $C_1$ and $C_2$ (see Figure 1),
identified with the initial curve $C$.

{\tt    \setlength{\unitlength}{0.92pt}}
\begin{figure} \label{constr}
   \epsfig{file=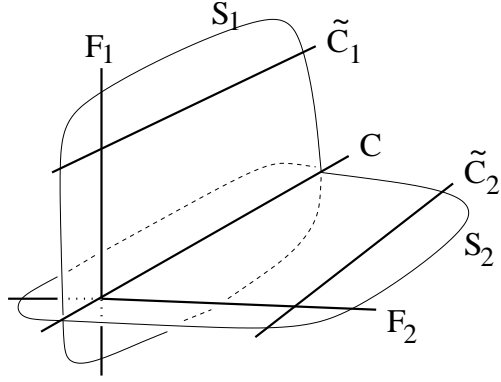, height= 50 mm}
   \caption{Construction of the 3-fold $M$}
\end{figure}

By construction, one has the following identification of the algebraic
normal bundles of $S_1$ and $S_2$ inside $M$:
\begin{equation} \label{identnorm}
    N_{S_i | M} \simeq H_i
\end{equation}

Combining this with the relations (\ref{secline}) we get:
\begin{equation} \label{intM}
  \left\{ \begin{array}{l}
     C \cdot_{M} S_i = \mbox{deg}_{C} N_{S_i | M}  = 
            \mbox{deg}_{C_i}H_i = -d_j\\
     F_i \cdot_{M} S_i = \mbox{deg}_{F_i} N_{S_i | M}  = 
            \mbox{deg}_{F_i}H_i = -x_i\\
     F_j \cdot_{M} S_i = \mbox{deg}_{F_j} N_{S_i | M}  = 
            \mbox{deg}_{F_j}H_i = 1
           \end{array} \right .
\end{equation}

In order to apply the criterion \ref{amplebij}, we want to find under
which conditions on the numbers $(d_1, d_2, x_1, 
x_2) \in \mathbf{Z}_{>0} ^4$, there exist pairs
$(\alpha_1, \alpha_2)\in \Z_{>0}^2$ such that the line bundle
$\mathcal{O}_M(-(\alpha_1S_1 + \alpha_2 S_2))$ is ample in
restriction to $S_1 \cup S_2$ (remember that we have already
imposed the restrictions (\ref{neg}) and (\ref{sign}) ).

By Kleiman's ampleness criterion \ref{klei} and Proposition 
\ref{amplecor}, this is equivalent to:
$$
  \left\{ \begin{array}{l}
       C \cdot_M (\alpha_1 S_1 + \alpha_2 S_2)<0\\
       F_1 \cdot_M (\alpha_1 S_1 + \alpha_2 S_2)<0\\
       F_2 \cdot_M (\alpha_1 S_1 + \alpha_2 S_2)<0
          \end{array} \right .
  \stackrel{(\ref{intM})}{\Longleftrightarrow} 
    \left\{ \begin{array}{l}
        \alpha_1 d_2 + \alpha_2 d_1 >0\\
        \alpha_1 x_1 -\alpha_2 >0\\
        \alpha_2 x_2 - \alpha_1 >0
             \end{array}
    \right.
$$

This in turn is equivalent to:

\begin{equation} \label{ineqampl}
  \frac{1}{x_1} < \frac{\alpha_1}{\alpha_2} < x_2
\end{equation}

The inequalities (\ref{ineqampl}) have solutions
$(\alpha_1, \alpha_2) \in \Z_{>0}^2$ if and only if at least one of
the numbers $x_1, x_2$ is $\geq 2$. It has solutions in both 
half-planes $\alpha_1 \geq \alpha_2$ and $\alpha_2 \geq \alpha_1$ if and only
if we have simultaneously $x_1\geq 2, \: x_2 \geq 2$. 

Combining this with Theorem \ref{amplecor} and Corollary
\ref{amplebij}, we get: 

\begin{proposition} \label{twocomp}
  Suppose that $\mathrm{deg}_C L_1 <0, \: \mathrm{deg}_C L_2
  <0$. Consider $x_1, x_2 \in 
  \mathbf{Z}_{>0}$. 
  If $x_1, x_2 \geq 1$ and at least one of them is $\geq 2$, then
  $S_1 \cup S_2$ is exceptional 
  in $M$.  Let then $\pi: (M, S_1 \cup S_2) \rightarrow (X,0)$ be
  the morphism which collapses $S_1 \cup 
  S_2$ in $M$. If moreover both $x_1$ and $x_2$ are $\geq 2$,
  then $S_1$ and $S_2$ are both essential components over
  $(X,0)$ and the Nash map $\mathcal{N}$ is bijective.
\end{proposition}

\begin{remark}
  We do not say that the Nash map is not bijective when one of the
  numbers $x_1, x_2$ is equal to $1$. But the
  method of the present paper does not allow us to decide it in
  general. 
\end{remark}

Our construction shows that the analytic germ $(X,0)$ defined in
Proposition \ref{twocomp} is uniquely determined by the choice of the
curve $C$, the line bundles $L_1, L_2$ and the numbers $x_1,
x_2$. That is why, when we want to recall these ingredients, we denote
it by $$(X_{C, L_1, L_2, x_1, x_2},0).$$
In the same way, we denote by $$M_{C, L_1, L_2, x_1, x_2}$$
the smooth algebraic manifold used to construct it.

\bigskip

\section{Analytic invariants of our families of examples}
\label{aninv}

In the introduction, we gave the list of the known examples of
higher-dimensional normal germs with isolated singularity which have a
bijective Nash map. It is natural to ask if the examples constructed
in the previous section are new, or cover partially the known ones. As
the normal quasi-ordinary germs are isomorphic to germs of simplicial
toric varieties (a result proved by the second author \cite[Theorem
5.1]{PP 04}, generalizing like this the hypersurface case treated by
Gonz{\'a}lez P{\'e}rez), and as   
in all our examples there are exactly two essential divisors, this
amounts to ask if some of them are analytically isomorphic to germs of
normal toric varieties. Through the propositions \ref{differ} and
\ref{toric}, we show that this is the case only when the curve $C$ is
rational.  
\medskip

The next proposition is a direct generalization of a result proved by
Nash \cite[page 35]{N 95}. It compares from the viewpoint of
birational algebraic geometry the centers on different
resolutions of a given essential divisor over $(X,0)$. 

\begin{proposition} \label{equicent}
 Let $(X,0)$ be a germ of normal analytic variety of dimension $n \geq
 2$, with isolated singularity. If $\pi_k : \tilde{X}_k \rightarrow X, \: k
 =1,2$ are two resolutions of $X$ and $E_k\subset \tilde{X}_k$ are  
 essential components corresponding to the same essential divisor
 over $(X,0)$, then $E_1 \times \mathbf{P}^{n-c_1 -1}$
 is birationally equivalent to $E_2 \times \mathbf{P}^{n-c_2 -1}$, where
 $c_k:= codim_{\tilde{X}_k}E_k$. 
\end{proposition}

\textbf{Proof:} 
  Denote by $\nu$ the essential divisor whose center on
  $\tilde{X}_k$ is $E_k$, for $k \in \{1,2\}$. Consider the morphism
  $\beta_k: B_{E_k}(\tilde{X}_k) \rightarrow \tilde{X}_k$ obtained by
  blowing-up $E_k$ in $\tilde{X}_k$, and the strict transform $D_k$
  of $E_k$ in $B_{E_k}(\tilde{X}_k)$. Then, $E_k$ is the center of
  $\nu$ in $B_{E_k}(\tilde{X}_k)$. It is birationally equivalent to
  $E_k \times \mathbf{P}^{n-c_k -1}$. Indeed, there exists a smooth
  Zariski open set $U_k\subset E_k$ whose normal bundle in
  $\tilde{X}_k$ is algebraically isomorphic to $U_k \times \mathbf{A}^{n -
  c_k}$, which shows that $\beta_k^{-1}(U_k) \simeq U_k \times
  \mathbf{P}^{n-c_k -1}$. But $D_k$ is birationally equivalent to
  $\beta_k^{-1}(U_k)$.

  Now consider the bimeromorphic map $\rho:= (\pi \circ \beta_2)^{-1}
  \circ (\pi \circ \beta_1): B_{E_1}(\tilde{X}_1) \rightarrow
  B_{E_2}(\tilde{X}_2)$. As the center of the valuation $\nu$ on
  $B_{E_k}(\tilde{X}_k)$ is the irreducible hypersurface $D_k$, this
  shows that the closure of $\rho(D_1)$ in $B_{E_2}(\tilde{X}_2)$ is
  equal to $D_2$. This means that $\rho$ realizes a birational
  equivalence between $D_1$ and $D_2$. The conclusion of the
  proposition follows.
\hfill $\Box$

\medskip

In order to analyze the germs $(X_{C, L_1, L_2, x_1, x_2},0)$
constructed in the previous section, we will use 
Proposition \ref{equicent} only through its Corollary
\ref{ruled}. Before stating it, let us introduce some notations.

Suppose that $(X,0)$ is a normal germ of 3-fold with isolated
  singularity. Consider a fixed resolution $\pi: \tilde{X} \rightarrow
  X$ of it. If $v$ is an 
  essential divisor over $(X,0)$, and $E_v$ is its center on
  $\tilde{X}$,  define its \emph{smooth representative} $R(v)$ to be:

$\bullet$ the unique 
  minimal model of $E_v$, if $E_v$ is a curve or a surface which is not
  birationally ruled;

$\bullet$ the curve $C$, if $E_v$ is birationally
  equivalent to $C \times \mathbf{P}^1$.

Recall from the introduction that $\mathcal{E}(X,0)$ denotes the set
of essential divisors over $(X,0)$. 

\begin{corollary} \label{ruled}
  The  collection $(R(v))_{v \in \mathcal{E}(X,0)}$ of abstract smooth 
  curves or minimal surfaces, parametrized by the set of essential
  divisors of $(X,0)$,  is independent of the choice of resolution. 
\end{corollary}

\textbf{Proof:} This is a direct consequence of the previous
proposition and of the fact that a non-birationally ruled surface has
a unique minimal model, whether if the smooth projective surfaces $C_1 \times
\mathbf{P}^1$ and $C_2 \times \mathbf{P}^1$ are birationally
equivalent, then the curves $C_1$ and $C_2$ are isomorphic  (see
B\u{a}descu \cite{B 01}).
\hfill $\Box$
\medskip

An immediate consequence of the corollary is:

\begin{proposition} \label{differ}
  1) If $(C, L_1, L_2, x_1, x_2)$ and $(C', L'_1, L'_2, x'_1,
     x'_2)$ are chosen such that $C, C'$ are non-isomorphic smooth
     projective curves, then the germs $(X_{C, L_1, L_2, x_1,
     x_2},0)$ and $(X_{C', L'_1, L'_2, x'_1, x'_2},0)$ are not
     analytically isomorphic.

  2) If $C$ is not rational, then $(X_{C, L_1, L_2, x_1,
     x_2},0)$ is not analytically isomorphic to a germ of toric variety.

\end{proposition}

\textbf{Proof:} 1) The set of smooth representatives of the essential
divisors of the germ $(X_{C, L_1, L_2, x_1,   x_2},0)$, each
representative being counted with its multiplicity, is equal to
$2C$. Point 1) follows then from Corollary \ref{ruled}. 

2) Any germ $(X,0)$ of toric variety has toric resolutions. The
irreducible components of the exceptional locus of such a resolution
are orbit closures, and in particular are rational varieties. This
shows that, when $X$ has dimension 3, 
all the smooth representatives of the essential divisors are rational
curves. Point 2) follows immediately. 
\hfill $\Box$
\medskip

\begin{proposition} \label{toric}

If $C$ is rational curve, then $(X_{C, L_1, L_2, x_1, x_2},0)$ is
analytically isomorphic to the germ of an affine normal toric variety
at the unique $0$-dimensional orbit.
\end{proposition}

\textbf{Proof:} If $C$ is rational, it is isomorphic to a toric curve,
and $S_1, S_2$ are isomorphic to toric surfaces, because the
only geometrically ruled surfaces over $\mathbf{P}^1$ are the
Hirzebruch surfaces. By recalling the shapes of the fans which define
the Hirzebruch surfaces and the way one gets the fan defining an orbit
closure from a given fan (see Fulton \cite{F 93}),  one sees that a
candidate fan 
$\Delta$ for a toric 3-fold isomorphic to $M$ and such that $S_1,
S_2, C$ are orbit closures should be as in  Figure 2.

{\tt    \setlength{\unitlength}{0.92pt}}
\begin{figure} \label{transv}
   \epsfig{file=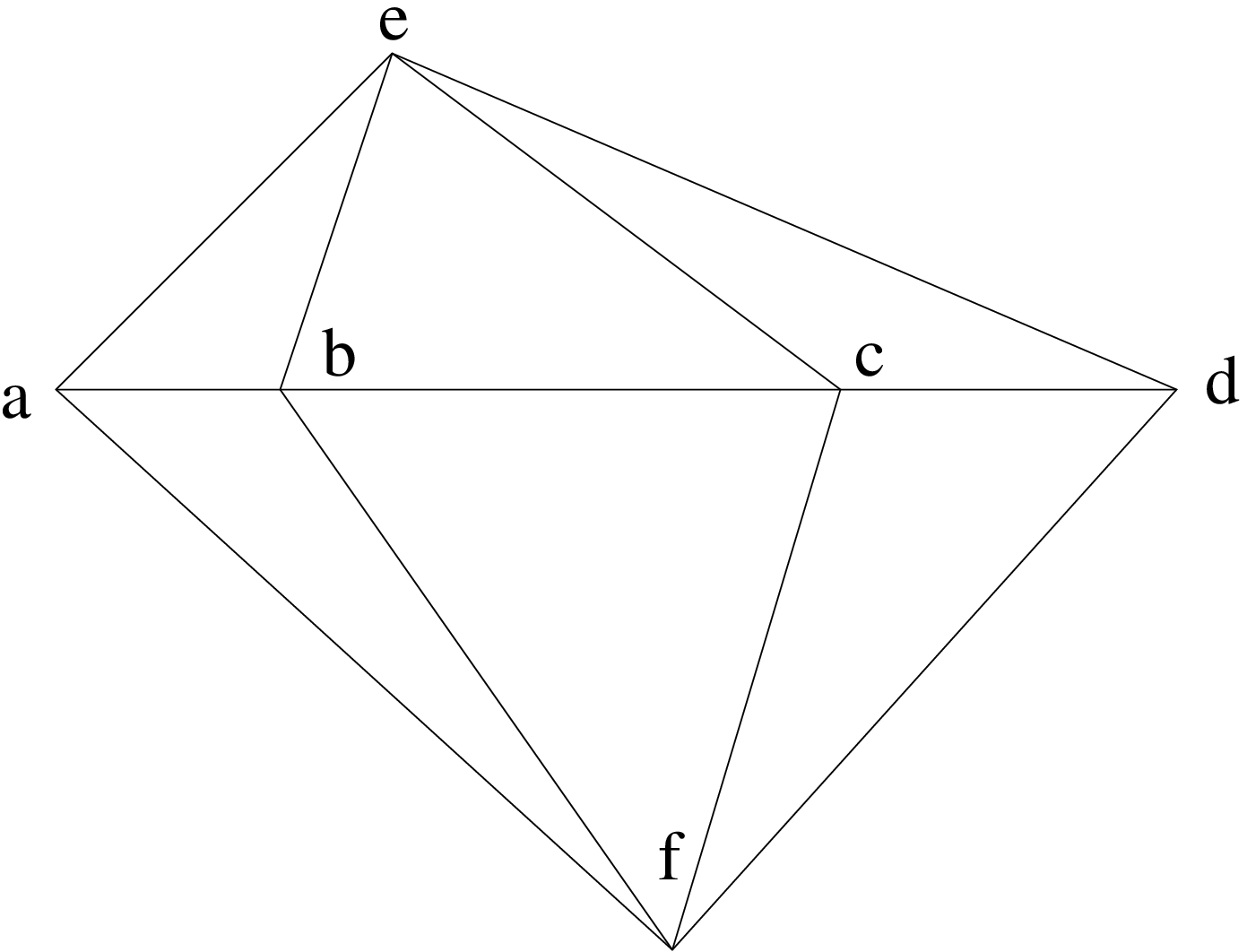, height= 50 mm}
   \caption{Figure illustrating the proof of Proposition \ref{toric}}
\end{figure}


In it, we have represented the intersections of the edges
$a,b,c,d,e,f$ of the fan (that is, its 1-dimensional cones) with a
transversal plane. The fan $\Delta$ lives inside the 3-dimensional
real vector space $N_{\R}$, where $N$ is the associated weight
lattice. For each edge $l$, we denote by $v_l$ the unique primitive
vector of $N$ situated on $l$. For each cone $\sigma$ of $\Delta$, we
denote by $O_{\sigma}$ the associated orbit and by $V_{\sigma}:=
\overline{O}_{\sigma}$ the orbit closure inside the toric variety
$\mathcal{Z}(N, \Delta)$. If $\sigma$ is strictly convex with edges
$l_1,...,l_n$, we write also $\sigma= \langle l_1,...,l_n\rangle$. 

As we want $\mathcal{Z}(N, \Delta)$ to be smooth, $\Delta$ must be a
regular fan. Moreover, we would like to get $S_1= V_b,
S_2=V_c, C= V_{\langle b, c \rangle}$ verifying the numerical
constraints (\ref{intM}). Those equations are equivalent in our toric
context with:

\begin{equation} \label{intor}
  \left\{ \begin{array}{l}
              V_{\langle b, c \rangle} \cdot V_b = - d_2\\
              V_{\langle b, c \rangle} \cdot V_c= - d_1\\
              V_{\langle b, e \rangle} \cdot V_b= - x_1\\
              V_{\langle c, e \rangle} \cdot V_c= - x_2
          \end{array} \right.
\end{equation}
where the intersection numbers are taken inside $\mathcal{Z}(N,
\Delta)$. The equalities \linebreak 
  $V_{\langle b, e \rangle} \cdot V_c=
V_{\langle c, e \rangle} \cdot V_b = 1$ are automatically satisfied,
as $\mathcal{Z}(N,\Delta)$ is smooth. 

We express the vectors $v_c, v_d, v_f$ in the basis $(v_a, v_b, v_e)$
of $N$. As we want $V_b, V_c$ to be Hirzebruch surfaces such that
$V_{\langle b, e \rangle}= V_b \cap V_c$ has negative
self-intersection in both of them, we require that $v_a, v_b, v_c,
v_d$ be coplanar. It is the matter of a simple computation to see that
this condition, combined with the requirement that $\Delta$ be
regular, shows the existence of $\alpha_1,..., \alpha_4 \in \Z$ such
that:
\begin{equation} \label{allvect}
  \left\{ \begin{array}{l}
              v_c = - v_a + \alpha_1 v_b\\
              v_d = - \alpha_2 v_a + (\alpha_1 \alpha_2 -1) v_b\\
              v_f = \alpha_3 v_a + \alpha_4 v_b - v_e
          \end{array} \right.
\end{equation}

In order to compute the intersection numbers of the left-hand side of
(\ref{intor}) from relations (\ref{allvect}), we use the general
formula (see Fulton \cite[Section 3.3]{F 93}): 
\begin{equation} \label{divtor}
  \mathrm{div}(\chi^m) = \sum_{l \in \Delta^{(1)}} (m, v_l) V_l, \: \forall \:
  m \in M
\end{equation}
where $\chi^m$ is the monomial corresponding to $m \in M$ (a character
of the associated algebraic torus) seen as a
rational function on $\mathcal{Z}(N, \Delta)$, and $\Delta^{(1)}$ is
the set of edges of $\Delta$. Here $M:= \mathrm{Hom}(N, \Z)$ denotes the
lattice of exponents of monomials. 

In our case, if $(v_a^*, v_b^*, v_c^*)$ denotes the basis of $M$ dual
of $(v_a, v_b, v_c)$, we get the following formulae by applying
(\ref{divtor}) to $m \in \{v_a^*, v_b^*\}$:
 $$  \left\{ \begin{array}{l}
               \mathrm{div}(\chi^{v_a^*})= V_a - V_c - \alpha_2 V_d + \alpha_3 V_f\\
               \mathrm{div}(\chi^{v_b^*})= V_b + \alpha_1 V_c + (\alpha_1 \alpha_2
               -1) V_d +\alpha_4 V_f 
          \end{array}  \right.$$
This implies:
 $$ \left\{ \begin{array}{l}
               0 = V_{\langle b, c \rangle}\cdot \mathrm{div}(\chi^{v_a^*})=
               -V_{\langle b, c \rangle}\cdot  V_c + \alpha_3\\
               0 = V_{\langle b, c \rangle}\cdot \mathrm{div}(\chi^{v_b^*})=
               V_{\langle b, c \rangle}\cdot  V_b + \alpha_1
               V_{\langle b, c \rangle}\cdot  V_c  + \alpha_4\\ 
               0 = V_{\langle b, e \rangle}\cdot \mathrm{div}(\chi^{v_b^*})=
               V_{\langle b, e \rangle}\cdot  V_b + \alpha_1\\
               0 = V_{\langle c, e \rangle}\cdot \mathrm{div}(\chi^{v_a^*})=
               -V_{\langle c, e \rangle}\cdot  V_c - \alpha_2
          \end{array}  \right.$$
By combining this with the equalities (\ref{intor}), we get:
  $$ \alpha_1= x_1, \: \alpha_2 =x_2, \: \alpha_3= -d_1, \: \alpha_4=
  d_2 +d_1 x_1$$
We have obtained the required decomposition of the vectors $v_c, v_d,
  v_f$ in the basis $(v_a, v_b, v_c)$ of $N$:
\begin{equation} \label{finform}
  \left\{ \begin{array}{l}
               v_c= -v_a + x_1 v_b\\
               v_d = -x_2 v_a + (x_1 x_2 -1) v_b\\
               v_f = -d_1 v_a + ( d_2 + d_1 x_1) v_b - v_e
          \end{array}  \right.
\end{equation}

We want to see now if this fan is a subdivision of a strictly convex
cone $\gamma$ in $N_{\R}$. This is equivalent to the fact that $a,
e,d,f$ are the edges of a strictly convex cone. After some routine
computations, one sees that the only non-trivial requirement is that a
linear form on $N_{\R}$ which vanishes on $v_d$ and $v_e$ takes
non-vanishing values of the same sign on $v_a$ and $v_f$. Such a
linear form is $m := x_2 v_b^* + (x_1 x_2 -1) v_a^* \: \in \:
M$. Then:
$$ \left\{ \begin{array}{l}
              (m, v_a) =x_1 x_2 -1\\
              (m,v_f) = d_1 + x_2 d_2
           \end{array} \right. $$
As $d_1 + x_2 d_2 >0$, we have to require that $x_1x_2 >1$. This is
precisely the condition we had found at the end of Section
\ref{examples}, ensuring that $S_1 \cup S_2$ is exceptional
in $M$ (see Proposition  \ref{twocomp}). We deduce that, if $x_1x_2
>1$, then $\gamma_{d_1, d_2, x_1, x_2} := \langle a,e,d,f\rangle$ is a
strictly convex cone whose edges are $a,e,d,f$.

Now, as $\mathcal{Z}(N, \Delta)$ is toric, we can easily show that it
is isomorphic to the manifold $M_{\mathbf{P}^1, L_1, L_2, x_1, x_2}$,
where $\mathrm{deg}_{\mathbf{P}^1}L_1 = -d_1, \:
\mathrm{deg}_{\mathbf{P}^1}L_2 = 
-d_2$. We deduce:
$$ (X_{\mathbf{P}^1, L_1, L_2, x_1, x_2},0) \simeq (\mathcal{Z}(N,
\gamma_{d_1, d_2, x_1, x_2}),0).$$ 
The proposition is proved.
\hfill $\Box$

\medskip

\end{document}